\documentclass[12pt,reqno,intlimits]{article}

\usepackage{amssymb,amsthm,amsmath,epsfig,graphicx}

\theoremstyle{remark}

\begin{document}
\begin{center} {\bf TWO-WEIGHTED INEQUALITIES FOR HARDY-LITTLEWOOD
MAXIMAL FUNCTIONS AND SINGULAR INTEGRALS IN $L^{p(\cdot)}$ SPACES}
\end{center}
\vskip+1cm

\begin{center}
{\bf Vakhtang Kokilashvili and Alexander Meskhi}
\end{center}

\vskip+1cm

{\bf Abstract.} Two--weight criteria of various type for the
Hardy--Littlewood maximal operator and singular integrals in
variable exponent Lebesgue spaces defined on the real line are
established. \vskip+0.5cm

\vskip+0.5cm


\vskip+0.5cm


\vskip+1cm

\vskip+0.5cm

\section*{Introduction}

We study the two-weight problem  for for Hardy--Littlewood maximal
functions and singular integrals in variable exponent Lebesgue
spaces $L^{p(\cdot)}$. In particular, we derive various type
two--weight criteria for the maximal functions and the Hilbert
transforms on the line. For a bounded interval   we assume that
the exponent $p$ satisfies the local log-H\"older continuity condition
and for the real line we require that $p$ is constant outside some
interval.  In the framework of variable exponent analysis such a
condition first appeared  in the paper \cite{Di1}, where the
author established the boundedness of the Hardy--Littlewood
maximal operator in $L^{p(\cdot)}({\Bbb{R}}^n)$. Unfortunately we
do not know whether the established criteria remain valid or not
when $p$ satisfies log--H\"older decay condition at infinity (see
\cite{CrFiNe} for this condition). It is known that the local log-H\"older continuity
condition  for the exponent $p$ together with the log-H\"older
decay condition guarantees the boundedness of operators  of
harmonic analysis in $L^{p(\cdot)}({\Bbb{R}}^n)$ spaces (see
\cite{CrFiNe}, \cite{Ne}, \cite{CaCrFi}, \cite{CrFiMaPe}).

The boundedness of the maximal, potential and singular operators
in $L^{p(\cdot)}({\Bbb{R}}^n)$ spaces was derived in the papers
\cite{Di1}, \cite{Di2}, \cite{DiRu}, \cite{CrFiNe}, \cite{Ne},
\cite{CrFiMaPe}, \cite{CaCrFi}. Weighted inequalities for
classical operators in $L^{p(\cdot)}_w$ spaces, were $w$ is a
power--type weight, were established  in the papers
\cite{KoSa1}-\cite{KoSa5}, \cite{SaVa}, \cite{SaSaVa}, \cite{DiSa}
etc,  while the same problems with general weights for Hardy,
maximal and fractional integral operators   were studied in
\cite{EdKoMeJFSA}-\cite{EdKoMe2}, \cite{KoMe4}, \cite{KoSa3},
\cite{Kop},  \cite{MaZe}, \cite{DiHa}. Moreover, in \cite{DiHa} a complete solution of the one--weight problem for maximal
functions defined on Euclidean spaces are given in terms of Muckenhoupt--type conditions. Finally we notice that in the paper \cite{EdKoMe2} modular--type sufficient conditions  governing the two--weight inequality for maximal and singular operators were established.

Throughout the paper $J$ denotes an interval (bounded or unbounded) in ${\Bbb{R}}$.

Let $p$ be a non--negative  function on ${\Bbb{R}}$. Suppose that
$E$ is a
 measurable subset of ${\Bbb{R}}$.   We use the following notation:
\begin{eqnarray*}
 p_-(E):= \inf_{E} p; \;\; p_+(E):= \sup_{E} p; \;\;\  p_-:= p_-({\Bbb{R}}); \;\; p_+:=
 p_+({\Bbb{R}}).
 \end{eqnarray*}

Assume that $1\leq p_-(J)\leq p_+(J)< \infty$. The variable
exponent Lebesgue space $L^{p(\cdot)}(J)$ (sometimes it is denoted
by $L^{p(x)}(J)$) is the class of all $\mu$-measurable functions
$f$ on $X$ for which $ S_p(f):= \int\limits_{J} |f(x)|^{p(x)} dx <
\infty. $ The norm in $L^{p(\cdot)}(J)$ is defined as follows:

$$ \|f\|_{L^{p(\cdot)}(J)} = \inf \{ \lambda>0: S_{p}(f/\lambda) \leq 1 \}. $$

It is known (see e.g. \cite{KoRa}, \cite{Sa1}, \cite{KoSa1}) that $L^{p(\cdot)}$  is a Banach space. For other
properties of $L^{p(\cdot)}$ spaces  we refer, e.g., to \cite{Sh},
\cite{KoRa}, \cite{Sa1}.

\vskip+0.1cm

Finally we point out that constants (often different constants in
the same series of inequalities) will generally be denoted by $c$
or $C$. The symbol $f(x) \approx g(x)$ means that there are
positive constants $c_1$ and $c_2$ independent of $x$ such that
the inequality $ f(x) \leq c_1 g(x) \leq c_2 f(x)$ holds.
Throughout the paper by the symbol $p'(x)$ is denoted the function
$p(x)/ (p(x)-1)$.

\section{Sawyer-type Condition for Maximal Operators in
$L^{p(x)}$ Spaces.}

\subsection{The case of bounded interval}

Let $J$ be bounded interval in $\mathbb{R}$ and let
$$(M^{(J)}_{\alpha}f)(x)=\sup\limits_{\substack{I \ni x \\  I\subset J}} \frac{1}{|I|^{1-\alpha}}\int\limits_{I}|f(y)|dy,
\;\;\; x\in J,$$ where $x\in J$ and $\alpha$ is a constant
satisfying the condition $0\leq \alpha <1$.

For a weight function $u$ we denote
$$ u(E):= \int\limits_E u(x) dx. $$
\vskip+0.1cm

{\bf Definition 1.1.}  Let $J$ be a bounded interval in
${\Bbb{R}}$. We say that a non--negative function $u$ satisfies
the doubling condition on $J$ ($u\in DC(J)$) if there is a
positive constant $b$ such that for all $x\in J$ and all $r$,
$0<r< |J|$,  the inequality
$$ u\big( I(x-2r, x+2r) \cap J\big) \leq b u\big( I(x-r, x+r)\cap J \big) $$ holds.

\vskip+0.1cm

{\bf Definition 1.2.} We say that  $p\in LH(J)$ ( $p$ satisfies
the local log-H\"older  condition)  if there is a positive
constant $c$ such that
$$  |p(x)-p(y)| \leq \frac{c}{-|x-y|} $$
for all $x,y \in J$ satisfying the condition $|x-y|\leq 1/2$.

\vskip+0.1cm

%Theorem 1
{\bf Theorem 1.1.} {\em  Let $1<p_{-}\leq p(x)\leq p_{+}<\infty$
and let the measure $d\nu(x)=w(x)^{-p'(x)}dx$ belongs to $DC(J).$
Suppose that $0\leq\alpha <1$ and that  $p\in LH(J)$ . Then the
inequality
$$\|v(\cdot) M^{(J)}_{\alpha}f\|_{L^{p(\cdot)}(J)}\leq c\|w(\cdot)f(\cdot)\|_{L^{p(\cdot)}(J)}$$  holds, if and only if there exist a positive constant  c such that for all interval $I$, $I\subset J$,
$$
    \int\limits_{I}(v(x))^{p(x)}(M^{(J)}_{\alpha}(w(\cdot)^{-p'(\cdot)}\chi_{I(\cdot)}))^{p(x)}dx\leq c\int\limits_{I}w^{-p'(x)}dx<\infty.
$$
} \vskip+0.2cm

To prove Theorem 1.1 we need some auxiliary statements.
\vskip+0.1cm

%Proposition 1
{\bf Proposition A. (\cite{SaWh}, Lemma 3.20)} {\em Let $s$ be a constant
satisfying the condition  $1<s<\infty$ and let $u\geq0$ on
$\mathbb{R}$. Suppose that  $\{Q_{i}\}_{i\in A}$ is a countable
collection of dyadic intervals in $\mathbb{R}$ and that
$\{a_{i}\}_{i\in A}$,$\{b_{i}\}_{i\in A}$ are sequences of
positive numbers satisfying the conditions:

$(i)\int\limits_{Q_{i}}u\leq a_{i} \;\; \text{for all}\; i\in A$;

$(ii)\sum\limits_{\{j\in A: Q_{j}\subset Q_{i}\}}b_{j}\leq ca_{i}
\;\; \text{for all} \; i\in A$.

Then there is a positive constant $c_s$ depended on $s$ such that
the inequality
\begin{equation*}
    \Bigg(\sum\limits_{i\in A}b_{i}\bigg(\frac{1}{a_{i}}\int\limits_{Q_{i}}gu\bigg)^{s}\Bigg)^{1/s}\leq c_{s}
    \bigg(\int\limits_{\mathbb{R}}g^{s}u\bigg)^{1/s}
\end{equation*}
holds for all non-negative functions $g$.}

\vskip+0.1cm

%Corollary 1
{\bf  Corollary A.} {\em Let $1<s<\infty$ and let $u$ be a non-negative measurable function on $\mathbb{R}$. Suppose that
$\{Q_{i}\}_{i\in A}$ is a a sequence  of dyadic cubes in
$\mathbb{R}^{n}$ and that $\{b_{i}\}_{i\in A}$ is a sequence of
positive numbers  satisfying the condition
$$\sum\limits_{\{j\in A: Q_{j}\subset Q_{i}\}}b_{j}\leq c u(Q_i).$$
Then there is a positive constant $c$ such that for all
non-negative functions $g$ the inequality
\begin{equation*}
    \sum\limits_{i\in A}b_{i}\bigg( \frac{1}{u(Q_i)}\int\limits_{Q_{i}} g u\bigg)^{s}\leq
    c\bigg(\int\limits_{\mathbb{R}}g^{s}u\bigg)^{1/s}
\end{equation*}
holds. }
 \vskip+0.1cm

{\bf Lemma A.} {\em Let $J$ be a bounded interval and let $1\leq  r_-(J)\leq r_+(J)<\infty$. Suppose that $r\in LH(J)$ and that the measure $\mu$ satisfies the condition $\mu \in DC(J)$. Then there is a positive constant
$c$ such that for all $f$, $\|f\|_{L^{r(\cdot)}(J, \mu)}\leq 1$, intervals
$I\subseteq J$ and $x\in I$ the inequality
$$ \bigg( \frac{1}{\mu(I)} \int\limits_I |f(y)| d\mu(y)  \bigg)^{r(x)}
\leq c \bigg[ \bigg( \frac{1}{\mu(I)} \int\limits_I
|f(y)|^{r(y)} d\mu(y)\bigg) +1\bigg]$$
holds.}

\vskip+0.1cm

{\em Proof.} We follow the idea of L. Diening \cite{Di1} (see also \cite{HaHaPe} for the similar statement in the case of metric measure spaces with doubling measure). We give the proof for completeness.

First recall that (see, e.g., \cite{HaHaPe}) since $J$ with the Euclidean distance and the measure $\mu$ is a bounded doubling space with the finite measure $\mu$ the condition $r\in LH(J)$ implies the following inequality:
$$ \big(\mu(I)\big)^{r_-(I)- r_+(I)} \leq C \eqno{(1.1)}$$
for all subintervals $I$ of $J$.

%Let  $x\in I$. Suppose that $\mu (I) >1/2$. Then
%$$ \bigg( \frac{1}{\mu(I)} \int\limits_I |f| d\mu(y)  \bigg)^{r(x)} \leq \bigg( \frac{1}{\mu (I)} \int\limits_I \big[ (f(y) \big)^{p(y)}+1\Big]
%d\mu(y) \bigg)^{p(x)} $$
%$$ c\bigg[ \bigg( \frac{1}{\nu (B)} \int\limits_B \Big( (f(y))^{p(y)}d\nu(y)
%+1\bigg)^{p(x)}+1\bigg] \leq  c \bigg[ \Big(
%f^{p(\cdot)}(\cdot)\Big)_{\nu, B}+1\bigg]$$

Assume  that $\nu B\leq 1/2$. By H\"older's inequality we have that

$$  \bigg( \frac{1}{\mu(I)} \int\limits_I |f(y)| d\mu(y)  \bigg)^{r(x)} \leq \bigg( \frac{1}{\mu (I) } \int_{I} |f(y)|^{r_-(I)} d\mu(y) \bigg)^{r(x)/r_-(I)} $$

$$ \leq c \mu(I)^{-r(x)/ r_-(I)} \bigg[\frac{1}{2}
\int\limits_{I}|f(y)|^{r(y)} d\mu(x) +  \frac{1}{2} \mu(I)\bigg]^{r(x)/r_-(I)}. $$

Observe now that the expression in brackets is less than or equal to $1$. Consequently, by (1.1) we find that

$$ \bigg( \frac{1}{\mu(I)} \int\limits_I |f(y)| d\mu(y)  \bigg)^{r(x)} \leq c   \mu(I)^{1- r(x) / r_-(I)} \bigg( \frac{1}{\mu(I)} \int\limits_I |f(y)|^{r(y)} d\mu(y)  +1 \bigg) $$

$$ \leq c   \mu(I)^{(r_-(I)- r_+(I))/r_-(I)} \bigg( \frac{1}{\mu(I)} \int\limits_I |f(y)|^{r(y)} d\mu(y)  +1 \bigg)
\leq c \bigg( \frac{1}{\mu(I)} \int\limits_I |f(y)|^{r(y)} d\mu(y)  +1 \bigg). $$

 The case $\mu(I)>1/2$ is trivial. $\Box$
\vskip+0.2cm

Suppose that $S$ is an interval in $\mathbb{R}$ and let us
introduce the dyadic maximal operator
 \begin{equation*}
    \big(M^{(d), S}_{\alpha}\big)f(x)=\sup\limits_{\substack{x\in I \\ I\in
    D(S)}}|I|^{\alpha-1}\int\limits_{I}|f(y)|dy,
 \end{equation*}
 where $0\leq \alpha <1$ and $D(S)$ is a dyadic lattice in $S$.

\vskip+0.1cm

To prove Theorem 1.1 we need the following statement:
 \vskip+0.1cm
%proof of which is given in \cite{Di1} for Euclidean spaces and in \cite{HaHaPe} for more general %case (see also \cite{KoMeAJM}):
% \vskip+0.1cm

% {\bf Theorem B.}  {\em Let $J$ be an interval (bounded or unbounded) and let $1< p_-(J)\leq p_+(J)<\infty$. Suppose that $\mu$ is a measure on $J$ %satisfying the doubling condition on $J$ $(\mu \in DC(J)$). Then $p\in LH(
%J)$ if and only if there is a positive constant $c$ such
%that for all intervals $I$, $\mu(I)\leq 1$, the inequality
%$$ \mu(I)^{p_-(I)- p_+(I)} \leq c $$
%holds.}
%\vskip+0.2cm

 %Lemma 1
{\bf Lemma 1.1.} {\em  Let $S$ be a bounded interval on $\mathbb{R}$ and
let $J$ be a subinterval of $S$. Suppose that
$\sigma(x):=w^{-p'(x)}$ belongs to the class $DC(J)$ and that
$p\in LH(J)$, where $1<p_{-}(J)\leq p(x)\leq p_{+}(J)<\infty$. Let
$0\leq \alpha <1$. If there is a positive constant $c$ such that
for all interval $I$, $I\subset J$,
  \begin{equation*}
    \int\limits_{I}(v(x))^{p(x)}\bigg(M^{(d),S}_{\alpha}\big(\chi_{I}(\cdot)\sigma(\cdot)\big)\bigg)^{p(x)}(x) dx\leq c\int\limits_{I}\sigma(x)dx<\infty,
  \end{equation*}
then the estimate
\begin{equation*}
\|v(\cdot)M^{(d),S}_{\alpha}\big(f(\cdot)\chi_{J}(\cdot)\big)\|_{L^{p(\cdot)}(J)}\leq
c\|w(\cdot)f(\cdot)\|_{L^{p(\cdot)}(J)}
\end{equation*}
holds. }

\vskip+0.1cm

{\em Proof.} Suppose that $\|f\|_{L^{p(\cdot)}_{w}(J)}\leq 1$. Assume that $f_1:=\chi_{J}f$. Let us introduce the set
\begin{equation*}
J_{k}=\{x\in S:2^{k}<(M^{(d),S}_{\alpha}f_1)(x)\leq 2^{k+1}\},\;\;
k\in\mathbb{Z}.
\end{equation*}
Suppose that for $k$, $J_k \neq \emptyset$, $\{I_{j}^k\}$ is a
maximal dyadic interval, $I_{j}^k\subset D(S)$, such that
$$
\frac{1}{|I_{j}^{k}|^{1-\alpha}}\int\limits_{I_{j}^{k}}|f_1(y)|dy>2^{k}.
\eqno{(1.2)}
$$
It is obvious that such  a maximal interval always exists. Now observe that

 $(i)\  \ \ \ \ \ \ \{I^{k}_{j}\}$ are disjoint for fixed $k$;

%$(ii)$
%$$|I^{k}_{j}|^{\alpha-1}\int\limits_{I^{k}_{j}}|f_1(y)|dy>2^{k}$$

$(ii)$
$$\overline{J}_{k}:= \{ x\in S: \big( M^{(d),S}_{\alpha} f_1\big)(x)> 2^{k} \}=  \cup_{j} I^{k}_{j}. $$
%where $\overline{I}_j^k:= I_j^k \cap S$.

\vskip+0.1cm

Indeed, $(i)$ holds because if $I_{i}^{k}\cap
I_{j}^{k}\not=\emptyset$, then $I_{i}^{k}\subset I_{j}^{k}$ or
$I_{j}^{k}\subset I_{i}^{k}$. Consequently, if $I_{i}^{k}\subset I_{j}^{k}$,
then $I_{j}^{k}$ is maximal interval for which (1.2) holds.

To see that $(ii)$ holds, observe that if $x\in
\overline{J}_{k}$, then $M^{(d),S}_{\alpha}f_1(x)\geq 2^{k}$. Hence,
there is a maximal dyadic interval $I_{j}^{k}$ containing  $x$
such that (1.2) hold for $I_{j}^{k}.$ Let now $x\in
\bigcup\limits_{j} I_{j}^{k}$.  Then $x\in I_{j_0}^{k}$
for some $j_0$. Hence,  $M^{(d),S}_{\alpha}f_1(x)> 2^{k}$ because
(1.2) holds for $I_{j_0}^{k}$.

Denote: 
\begin{equation*}
E_{j}^{k}:= I_{j}^{k}\backslash \{x\in
S:M^{(d),S}_{\alpha}f_1(x)> 2^{k+1}\}.
\end{equation*}
Then $E_{j}^{k}= I_{j}^{k}\cap J_{k}$. Indeed, if $x\in
E_{j}^{k}$, then $x\in I_{j}^{k}$ and  $M^{(d),S}_{\alpha}f_1(x) \leq 
2^{k+1}$.  Hence, by (1.2) we find that 
\begin{equation*}
    2^{k}<|I_{i}^{k}|^{\alpha -1}\int\limits_{I_{j}^{k}}|f_1(y)|dy\leq M^{(d),S}_{\alpha}f_1(x)\leq 2^{k+1}.
\end{equation*}
This means that $x\in  I_{j}^{k}\cap J_{k}$.
Let now  $x\in I_{j}^{k}\cap J_{k}$. Then  obviously 
$M^{(d),S}_{\alpha}f_1(x)\leq 2^{k+1}$. Consequently,  $x\in E_{j}^{k}.$\\
Observe that $\{E_{j}^{k}\}$ are disjoint for every $j,k$ because, as we have seen,
\begin{equation*}
E_{j}^{k}=\{x\in
I_{j}^{k}:2^{k}<M^{(d),S}_{\alpha}f_1(x)\leq 2^{k+1}\}.
\end{equation*}
Also, $ E_{j}^{k}\subset I_{j}^{k} $. Assume that
$\|w(\cdot)f_1(\cdot)\|_{L^{p(\cdot)}(S)}\leq 1$. Denote:
$$v_1:= v \chi_{J}, \;\;\;  \sigma_1:= \sigma \chi_{J}.$$
By the arguments observed above and using Lemma A with $r(\cdot)=p(\cdot)/p_-$ and the measure $d\mu(x)= \sigma(x) dx$  we have that
\begin{eqnarray*}
&&\int\limits_{J}(v(x))^{p(x)}\bigg(M^{(d),S}_{\alpha} f_1 \bigg)^{p(x)}(x)dx  \\&=&
\int\limits_{S}(v_1(x))^{p(x)}\bigg(M^{(d),S}_{\alpha} f_1 \bigg)^{p(x)}(x)dx\\
&\leq&\sum\limits_{j,k}\int\limits_{E_{j}^{k}}(v_1(x))^{p(x)}2^{(k+1)p(x)}dx\\
&\leq&
c\sum\limits_{j,k}\int\limits_{E_{j}^{k}}(v_1(x))^{p(x)}\bigg(\frac{1}{|I_{j}^{k}|^{1-\alpha}}
\int\limits_{I_{j}^{k}}|f_1(y)|dy\bigg)^{p(x)}dx\\
&=&c\sum\limits_{j,k}\int\limits_{E_{j}^{k}}(v_1(x))^{p(x)}
\bigg(\frac{\sigma(I_{j}^{k}\cap
J)}{|I_{j}^{k}|^{1-\alpha}}\bigg)^{p(x)}\bigg(\frac{1}{\sigma(I_{j}^{k}\cap
J)}\int\limits_{I_{j}^{k}}
\Big|\frac{f_1}{\sigma}\Big|\sigma \bigg)^{p(x)}dx\\
&=&c\sum\limits_{j,k}\int\limits_{E_{j}^{k}}(v_1(x))^{p(x)}
\bigg(\frac{\sigma(I_{j}^{k}\cap
J)}{|I_{j}^{k}|^{1-\alpha}}\bigg)^{p(x)}\bigg(\frac{1}{\sigma(I_{j}^{k}\cap
J)}
\int\limits_{I_{j}^{k}}\Big|\frac{f_1}{\sigma}\Big|\sigma \bigg)^{p(x)}dx\\
&\leq&c\sum\limits_{j,k}\Bigg(\int\limits_{E_{j}^{k}}(v_1(x))^{p(x)}
\bigg(\frac{\sigma(I_{j}^{k}\cap J)}{|I_{j}^{k}|^{1-\alpha}}
\bigg)^{p(x)}dx\Bigg)\bigg(\frac{1}{\sigma(I_{j}^{k}\cap
J)}\int\limits_{I_{j}^{k}}
\Big|\frac{f_1(y)}{\sigma(y)}\Big|^{\frac{p(y)}{p_{-}}}\sigma (y) dy\bigg)^{p_-}\\
&+&c\sum\limits_{j,k}\Bigg(\int\limits_{E_{j}^{k}}(v_1(x))^{p(x)}
\bigg(\frac{\sigma(I_{j}^{k}\cap J)}{|I_{j}^{k}|^{1-\alpha}}
\bigg)^{p(x)}dx\Bigg)\\
&\equiv&
c\bigg(\sum\limits_{j,k}A_{j}^{k}+\sum\limits_{j,k}B_{j}^{k}\bigg).
\end{eqnarray*}
Notice that  the sign of sum is taken over all those $j$ ad $k$
for which $\sigma(I_{j}^{k}\cap J)>0).$

To use Corollary A observe that

\begin{eqnarray*}
% \nonumber to remove numbering (before each equation)
     && \sum\limits_{\substack{I_{j}^{k}\subset I_{i}\\I_j^k,I_i\in D(S)}}\int\limits_{E_{j}^{k}}(v_1(x))^{p(x)}
     \bigg(\frac{\sigma(I_{j}^{k}\cap J)}{|I_{j}^{k}|^{1-\alpha}}\bigg)^{p(x)}dx \\
        &\leq& \sum\limits_{I_{j}^k\subset I_{i}}\int\limits_{E_{j}^{k}}(v_1(x))^{p(x)}\bigg( M^{(d),S}_{\alpha}(\chi_{I_{i}\cap J}\sigma)\bigg)^{p(x)}(x)dx\\
&\leq&\int\limits_{I_{i}} (v_1(x))^{p(x)}  \bigg( M^{(d),S}_{\alpha}(\chi_{I_{i}\cap J}\sigma)\bigg)^{p(x)}(x)dx\\
&\leq &c\int\limits_{I_i\cap J}\sigma(x)dx = c \int\limits_{I_i}
\sigma_1(x) dx.
\end{eqnarray*}
Now Corollary A implies that
       \begin{eqnarray*}
       % \nonumber to remove numbering (before each equation)
         \sum\limits_{j,k}A_{j}^{k} &=& \sum\limits_{j,k} \bigg( \int\limits_{E_{j}^{k}}(v_1(x))^{p(x)}
         \bigg(\frac{\sigma(I_{j}^{k}\cap J)}{|I_{j}^{k}|^{1-\alpha}}
         \bigg)^{p(x)}dx\Bigg) \bigg(\frac{1}{\sigma_1(I_{j}^{k})}\int\limits_{I_{j}^{k}}
         \Big|\frac{f_1(y)}{\sigma(y)}\Big|^{\frac{p(x)}{p_-}}\sigma_1(y) dy \bigg)^{p_-} \\
          &\leq& c \int\limits_{S}|f_1(x)|^{p(x)}\sigma(x)^{-p(x)}\sigma_1(x)dx =  c \int\limits_{S} |f_1(x)|^{p(x)} w^{p(x)} dx \leq c.
       \end{eqnarray*}
   For the second term we have that
   \begin{eqnarray*}
   % \nonumber to remove numbering (before each equation)
    \sum\limits_{j,k}B_{j}^{k} &=& \sum\limits_{j,k}\int\limits_{E_{j}^{k}}(v_1(x))^{p(x)}\bigg(\frac{\sigma(I_{j}^{k}\cap J)}{|I_{j}^{k}|^{1-\alpha}}\bigg)^{p(x)}dx \\
      &\leq& \sum\limits_{j,k}\int\limits_{E_{j}^{k}}(v_1(x))^{p(x)}\bigg(M^{(d),S}_{\alpha}(\chi_{J}\sigma)\bigg)^{p(x)}(x)dx\\
      &=&\int\limits_{J}(v(x))^{p(x)}\bigg(M^{(d),S}_{\alpha}(\chi_{J}\sigma)\bigg)^{p(x)}(x)dx\\
      &\leq& c\int\limits_{J}\sigma(x)dx<\infty.\\
   \end{eqnarray*}
   Finally we conclude that
   \begin{equation*}
    \|v(\cdot)\big(M^{(d),S}_{\alpha} f_1 \big)(\cdot)\|_{L^{p(\cdot)}(J)}\leq c
   \end{equation*}
   for $\|w(\cdot)f(\cdot)\|_{L^{p(\cdot)}(J)}\leq 1.$ $\;\;\; \Box$

 \vskip+0.2cm

{\em Proof of Theorem} 1.1.  {\em Sufficiency}. Let us take an
interval $S$ containing $J.$  Without loss of generality we can
assume that $S$ is a maximal dyadic interval and that $|J|\leq \frac{|S|}{8}$. Further, suppose also that 
$J$ and $S$ have one and the same center. Without loss of generality assume that
$|S|=2^{m_0}$ for some integer $m_0$. Then every interval
$I\subset J$ has the length $|I|$ less than or equal to
$2^{m_0-3}$ . Assume that $|I|\in[2^j,2^{j+1})$ for some $j$,  $j\leq
m_0-4$. Let us introduce the set
\begin{equation*}
F=\{t\in (-2^{m_0-4},2^{m_0-4}):\; \text{there is} \;\; I_1 \in
D(S)-t,I\subset I_1\subset S, |I_1|=2^{j+1}\}.
\end{equation*}
The simple geometric observation (see also  \cite{GaRu}, p. 431) shows
that $|F|\geq 2^{m_0-4}$.

Further, let 
  \begin{equation*}
    (K_tf)(x):=\sup\limits_{\substack{ S\supset I_1 \ni x \\I_1\in D(S)-t}}\frac{1}{|I_1|^{1-\alpha}}\int\limits_{I_1}|f_1|,\;\; t\in F,
  \end{equation*}
where $f_1=\chi_{J}f$. Then
for $x$ ($x\in J$) there  exist  $I\ni x$, $I\subset J$  such that
$$ |I|^{\alpha-1}\int\limits_{I}|f_1|>\frac{1}{2} (M^{(J)}_{\alpha}f_1)(x).$$ For the interval $I$, we have that $|I|\in[2^j,2^{j+1})$, $j \leq m_0-4$. Therefore for $t\in F$, there is an interval $I_1$, $I_1\in D(S)-t, I\subset I_1 \subset S,|I_1|=2^{j+1}$, such that 
$$
|I|^{\alpha-1}\int\limits_{I}|f_1|\leq\frac{c}{|I_{1}|^{1-\alpha}}\int\limits_{I_1}|f_1|.
$$
 Hence,
\begin{equation*}
(M_{\alpha}^{(J)}f)(x)\leq c(K_t f_1)(x), \;\; \hbox{for
every}\;\;  t\in F, x\in J,
\end{equation*}
with the positive constant $c$ depending only on $\alpha$.
Consequently,
\begin{eqnarray*}
% \nonumber to remove numbering (before each equation)
(M^{(J)}_{\alpha}f)(x) &\leq&  \frac{1}{|F|}\int\limits_{F}(K_tf_1)(x)dt \\
&\leq& \frac{c}{|I(0,2^{m_0-4})|}
\int\limits_{I(0,2^{m_0-4})}(K_tf_1)(x)dt.
\end{eqnarray*}

Suppose that $\|w(\cdot)f(\cdot)\|_{L^{p(\cdot)}(J)}\leq 1$. Then
by Lemma 1.1 we have that
\begin{eqnarray*}
% \nonumber to remove numbering (before each equation)
S_t &:=&\int\limits_{J}(v(x))^{p(x)}\big((K_tf_1)(x)\big)^{p(x)}dx  \\
&=&\int\limits_{J}(v(x))^{p(x)}\bigg(\sup\limits_{\substack{S
\supset I_1\ni x\\I_1\in D(S)-t}}
\frac{1}{|I_{1}|}\int\limits_{I_1}|f_1|\bigg)^{p(x)}dx  \\
&=& \int\limits_{J+t}(v_t(x))^{p(x-t)}\bigg(\sup\limits_{\substack{S \supset I_1\ni x\\I_1\in D(S)}}|I_1|^{\alpha-1}\int\limits_{I_1}\chi_{J}(s-t)f_1(s-t)ds\bigg)^{p(x-t)}dx\\
&=&\int\limits_{J+t}(v_t(x))^{p_t(x)}\bigg(\sup\limits_{\substack{I_1\ni x\\I_1\in D(S)}}|I_1|^{\alpha-1}\int\limits_{I_1}\chi_{J+t}(s)f_1(s-t)ds\bigg)^{p_t(x)}dx\\
&=& \int\limits_{J+t}(v_t(x))^{p_t(x)}\bigg(M^{(d),S}_{\alpha}\big(\chi_{J+t}(\cdot)f_1(\cdot-t)\big)\bigg)^{p_t(x)}dx\\
&\leq&c
\end{eqnarray*}
provided that
\begin{equation*}
\int\limits_{J+t}(w_t(x))^{p_t(x)}(f_1(x-t))^{p_t(x)}dx=\int\limits_{J}w(x)|f(x)|^{p(x)}dx\leq
1,
\end{equation*}
where $v_t(x)=v(x-t)$, $w_t(x)= w(x-t)$, $p_t(x)= p(x-t)$. To
justify this conclusion we need to check  that for every $I,I
\subset J+t$,
$$ \int\limits_{I}(v_t(x))^{p_t(x)}\bigg( M^{(d),S}_{\alpha}(\sigma_t\chi_I)(x)\bigg)^{p_t(x)}dx\leq c\int\limits_{I}\sigma_t(x)dx<\infty,
$$
where the positive constant $c$ is independent of $I$ and $t$. Indeed, observe that
\begin{eqnarray*}
% \nonumber to remove numbering (before each equation)
   &&\int\limits_{I}(v_t(x))^{p_t(x)}\bigg( M^{(d),S}
   _{\alpha}(\sigma_t\chi_I)(x)\bigg)^{p_t(x)}dx  \\
   &=& \int\limits_{I}(v_t(x))^{p_t(x)}\Bigg(\sup\limits
   _{\substack{ I_1\ni x\\I_1\in D(S)}}|I_1|^
   {\alpha-1}\int\limits_{I_1}\chi_I(s)
   \sigma(s-t)ds\Bigg)^{p_t(x)}dx \\
   &=& \int\limits_{I}(v_t(x))^{p_t(x)}
   \Bigg(\sup\limits_{\substack{ I_1-t
   \ni x-t\\I_1\in D(S)}}|I_1-t|^
   {\alpha-1}\int\limits_{I_1-t}\chi_I(s+t)\sigma(s)ds\Bigg)^{p_t(x)}dx \\
   &=& \int\limits_{I-t}(v(x))^{p(x)}
   \Bigg(\sup\limits_{\substack{ I_1
   \ni x\\I_1\in D(S)-t}}|I_1|^
   {\alpha-1}\int\limits_{I_1}
   \chi_{I-t}(s)\sigma(s)ds\Bigg)^{p(x)}dx  \\
   &\leq&\int\limits_{I-t}(v(x))^{p(x)}
   \bigg(M^{(J)}_{\alpha}(\chi_{I-t}\sigma)
   \bigg)^{p(x)}(x) dx\leq\int\limits_{I-t} \sigma(x)dx  \\
   &=& \int\limits_{I} \sigma_t(x) dx<\infty.
\end{eqnarray*}
Further, let  $g\in L^{p'(\cdot)}(J)$  with
$\|g\|_{L^{p'(\cdot)}(J)}\leq 1$.  Then we find that

\begin{eqnarray*}
% \nonumber to remove numbering (before each equation)
&&\int\limits_{J}(M^{(J)}_{\alpha}f)(x)v(x)g(x)dx \\
&\leq& \int\limits_{J}\bigg(\frac{1}{|I(0,2^{m_0-4})|}\int\limits_{I(0,2^{m_0-4})}(K_t f_1)(x)dt\bigg)v(x)g(x)dx \\
&\leq& \frac{1}{|I(0,2^{m_0-4})|}\int\limits_{I(0,2^{m_0-4})}\bigg( \int\limits_{J}(K_t f_1)(x)g(x)v(x)dx\bigg)dt \\
&\leq& \frac{1}{|I(0,2^{m_0-4})|}\int\limits_{I(0,2^{m_0-4})}\|(K_t f_1)v\|_{L^{p(\cdot)}(J)} \|g\|_{L^{p'(\cdot)}(J)}dt \\
&\leq&  c,
\end{eqnarray*}
provided that  $\|f\|_{L^{p(\cdot)}_{w}(J)}\leq 1$.

Finally we conclude that $ \|(M^{(J)}_{\alpha}f)
v\|_{L^{p(\cdot)}(J)} \leq c $ if $\|fw\|_{L^{p'(\cdot)}(J)}\leq
1.$

Sufficiency is proved.

{\em Necessity.} Let $f_I(t)=\chi_{I}(t)w^{-p'(t)}(t)$. Suppose
that $\beta=\|w^{-1}(\cdot)\|_{L^{p'(\cdot)}(J)}\leq1.$ We have
that
\begin{equation*}
    \big\|v(\cdot)(M^{(J)}_{\alpha}f)^{p(\cdot)}(\cdot)\big\|_{L^{p(\cdot)}(J)}\geq\big\|\chi_I(\cdot)
    v(\cdot)\big(M^{(J)}_{\alpha} \big(w^{-p'(\cdot)}(\cdot)\chi_{I}(\cdot)\big)\big)(\cdot)\big\|
    _{L^{p(\cdot)}(J)}=: A.
\end{equation*}
Hence, by the boundedness of $M_{\alpha}^{(J)}$, Lemma $B$ (recall that the measure $d\nu(x)=w(x)^{-p'(x)}dx$ satisfies the doubling condition) and the fact that $1/p \in LH(J)$ we find that

\begin{eqnarray*}
     % \nonumber to remove numbering (before each equation)
&A=&  \big \|\chi_I(\cdot)
v(\cdot) M^{(J)}_{\alpha}\big(w^{-p'(\cdot)}(\cdot)\chi_{I}(\cdot)\big)(\cdot) \big\|_{L^{p(\cdot)}(J)} \\
&\leq& c \big\|w(\cdot)w^{-p'(\cdot)}(\cdot)\chi_{I}(\cdot)\big\|_{L^{p(\cdot)}(J)} \\
&\leq&c\bigg(\int\limits_{I}w^{-p'(x)p(x)}(x)w^{p(x)}(x)dx\bigg)^{1/p_+(I)}  \\
%&\leq&c\bigg(\int\limits_{I}w^{(1-p'(x))p(x)}(x)dx\bigg)^{1/p_+(I)}  \\
%&\leq& c\bigg(\int\limits_{I}w^{-p'(x)}(x)dx\bigg)^{1/p_+(I)} \\
%&=&c\bigg(\int\limits_{I}w^{-p'(x)}(x)dx\bigg)^{(\frac{1}{p})_-(I)}  \\
&\leq& \bar{
c}\bigg(\int\limits_{I}w^{-p'(x)}(x)dx\bigg)^{\frac{1}{p_-(I)}}
\leq\bar{c}.
\end{eqnarray*}
On the other hand,
\begin{eqnarray*}
% \nonumber to remove numbering (before each equation)
&A=&\bar{c}\bigg\|\frac{1}{\bar{c}}\chi_I(\cdot)v(\cdot)
M^{(J)}_{\alpha}\big(w^{-p'(\cdot)}\chi_I(\cdot)\big)(\cdot)\bigg\|_{L^{p(\cdot)}(J)}   \\
&\geq&\bar{c}\bigg(\int\limits_{I}(\bar{c})^{-p(x)}(v(x))^{p(x)}
\bigg[M^{(J)}_{\alpha}\big(w^{-p'(\cdot)}\chi_I(\cdot)\big)\bigg](x)dx\bigg)^{\frac{1}{p_-(I)}}\\
&\geq&c\bigg[\int\limits_{I}(v(x))^{p(x)}\bigg(M^{(J)}_{\alpha}\big(w^{-p'(\cdot)}
\chi_I(\cdot)\big)(x)\bigg)^{p(x)}dx\bigg]^{\frac{1}{p_-(I)}}.
 \end{eqnarray*}

Summarizing these inequalities we conclude that

\begin{equation*}
    \int\limits_{I}(v(x))^{p(x)}\bigg(M^{(J)}_{\alpha}\big(w^{-p'(\cdot)}\chi_I(\cdot)\big)(x)\bigg)^{p(x)}dx\leq
    c\int\limits_{I}w^{-p'(x)}(x)dx<\infty.
\end{equation*}
Suppose now that $\beta \geq1$. Let us take
$$f(t)=\frac{w^{-p'(t)}(t)\chi_I(t)}{\beta}.$$ Then
$$\|f_I(\cdot)w(\cdot)\|_{L^{p(\cdot)}(J)}=\frac{\|w^{1-p'(\cdot)}(\cdot)\chi_I(\cdot)\|_{L^{p(\cdot)}(J)}}{\beta}\leq1.$$
Arguing as above we have desire result. It remains to show that
$$A:=\int\limits_{J}w^{-p'(x)}(x)dx<\infty.$$ Suppose that $A=\infty.$
Then $\|w^{-1}(\cdot)\|_{L^{p'(\cdot)}(J)}=\infty.$
Hence, there exist a function $g$, $\|g\|_{L^{p(\cdot)}(J)},g\geq
0$ such that
$$\int\limits_{J}g(x)w^{-1}(x)dx=\infty. $$
Let $f(x)=g(x)w^{-1}(x)$. Then
\begin{equation*}
    \bigg\|v(\cdot)\big(M^{(J)}_{\alpha}f\big)(\cdot)\bigg\|_{L^{p(\cdot)}(J)}
    \geq\bigg(\int\limits_{J}w^{-1}(x)g(x)\bigg)\bigg\|v(\cdot)
    |J|^{\alpha-1}\bigg\|_{L^{p(\cdot)}(J)}=\infty,
\end{equation*}
while
$$\|fw\|_{L^{p(\cdot)}(J)}=\|g\|_{L^{p(\cdot)}(J)}<\infty.$$ $\Box$
\vskip+0.2cm

%Corollary 2

{\bf Corollary 1.1.} {\em  Let $J$ be a bounded interval and let $1<p_-(J)\leq p(x)\leq p_+(J)<\infty$ and
let $0\leq \alpha<1$. Assume that $p\in LH(J)$ then the inequity
\begin{equation*}
    \big\|v(\cdot)\big(M^{(J)}_{\alpha}f\big)(\cdot)\big\|_{L^{p(\cdot)}(J)}\leq c \|f\|_{L^{p(\cdot)}(J)} \;\;( \text{Trace inequality})
\end{equation*}
holds if and only if
$$\sup\limits_{I,I\subset J}\frac{1}{|I|}
\int\limits_{I}(v(x))^{p(x)}|I|^{\alpha p(x)}dx<\infty. $$ }
\vskip+0.1cm

{\em Proof.} {\em Sufficiency.} By Theorem 1.1 it is enough to see
that

$$\big(M^{(J)}_{\alpha}\chi_{I}\big)(x)\leq |I|^{\alpha}  \ \ \ \ \hbox{for}\ \ \ \ \  x\in I.$$ This is true because of the following estimates:

$$\sup\limits_{\substack{S,S\subset J\\S\ni x}}|S|^{\alpha-1}\int\limits_{S}\chi_{I}\leq\sup\limits_{\substack{S\cap I\ni x\\S\subset J}}|S\cap I|^{\alpha-1}\int\limits_{S\cap I}dx=
\sup\limits_{\substack{S\cap I\ni x\\S\subset J}}|S\cap
I|^{\alpha}=|I|^{\alpha}. $$

{\em Necessity} follows by choosing the appropriate test functions in the trace inequality.
$\Box$

\subsection{The case of unbounded interval}

Now we derive criteria for the two--weight inequality for the following maximal operators:
\begin{equation*}
\bigg(M^{(\mathbb{R}_+)}_{\alpha}f\bigg)(x)=\sup\limits_{h>0}\frac{1}{h^{1-\alpha}}\int\limits_{(x-h,x+h)\cap\mathbb{R}_+}|f(y)|dy
\end{equation*}
and
\begin{equation*}
  \bigg(M^{(\mathbb{R})}_{\alpha}f\bigg)(x)=\sup\limits_{h>0}\frac{1}{h^{1-\alpha}}\int\limits_{x-h}^{x+h}|f(y)|dy,
\end{equation*}
where $0\leq  \alpha<1$.

In the sequel we will assume that $v^{p(\cdot)}(\cdot)$  and  $w^{-p'(\cdot)}(\cdot)$ are a.e. positive locally integrable function.
\vskip+0.2cm
{\bf Theorem 1.2.} {\em Let $0\leq \alpha<1$, $1<p_-({\Bbb{R}}_+)\leq p\leq p_+({\Bbb{R}}_+)<\infty$  and let $p\in LH(\mathbb{R}_+)$. Suppose that there is a bounded interval $[0,a]$ such that $w^{-p'(\cdot)}(\cdot)\in DC([0,a])$ and $p\equiv p_c\equiv $const outside $[0,a]$. Then the inequity
\begin{equation*}
    \|v M^{(\mathbb{R}_+)}_{\alpha}f\|_{L^{p(\cdot)}(\mathbb{R}_+)}\leq\|wf\|_{L^{p(\cdot)}(\mathbb{R}_+)},
\end{equation*}
holds  if and only if there is a positive constant $b$ such that for all bounded intervals $I\subset\mathbb{R}_+ $,
\begin{equation*}
    \|v M^{(\mathbb{R}_+)}_{\alpha}(w^{-p'(\cdot)}\chi_{I})\|_{L^{p(\cdot)}(I)}\leq c\|w^{1-p'(\cdot)}\|_{L^{p(\cdot)}(I)}<\infty.\eqno{(1.3)}
\end{equation*}}

\begin{proof} {\em Sufficiency.} Suppose that $\|wf\|_{L^{p(\cdot)}(\mathbb{R}_+)}<\infty$. We will show that $\|vM^{(\mathbb{R}_+)}_{\alpha}\|_{L^{p(\cdot)}(\mathbb{R}_+)}<\infty.$

Represent $ M^{(\mathbb{R}_+)}_{\alpha}f (x)$ as follows:
\begin{eqnarray*}
% \nonumber to remove numbering (before each equation)
  && M^{(\mathbb{R}_+)}_{\alpha}f(x) =\chi_{[0,a]}(x)M^{(\mathbb{R}_+)}_{\alpha} \big( f\cdot\chi_{[0,a]}\big) (x)\\
  &&
   + \chi_{[0,a]}(x)M^{(\mathbb{R}_+)}_{\alpha}\big(f\cdot\chi_{(a,\infty)}\big)(x)
   +\chi_{(a,\infty)}(x)M^{(\mathbb{R}_+)}_{\alpha}\big(f \cdot\chi_{[0,a]}\big)(x)\\
   &&+ \chi_{(a,\infty)}(x) M^{(\mathbb{R}_+)}_{\alpha}\big( f \cdot\chi_{(a,\infty)}\big)(x) \\
   &=:& M^{(1)}_{\alpha}f(x)+M^{(2)}_{\alpha}f(x)+M^{(3)}_{\alpha}f(x)
   + M^{(4)}_{\alpha}f(x).
\end{eqnarray*}
Since   $\|wf\|_{L^{p(\cdot)}(\mathbb{R}_+)}<\infty$ we have that  $\|wf\|_{L^{p(\cdot)}([0,a])}<\infty$. Applying now Theorem 1.1 we find that $\|vM^{(1)}_{\alpha}f\|_{L^{p(\cdot)}(\mathbb{R}_+)}<\infty$.
Further, observe that
\begin{equation*}
M^{(2)}_{\alpha}f(x)\leq\sup\limits_{h>a-x}\frac{1}{h}\int\limits_{a}^{x+h}|f(y)|dy
\leq\big(M^{(\mathbb{R}_+)}_{\alpha}f\big)(a) <\infty.
\end{equation*}
Hence,
\begin{equation*}
    \|vM^{(2)}_{\alpha}f\|_{L^{p(\cdot)}(\mathbb{R}_+)}\leq\big(M^{(\mathbb{R}_+)}_{\alpha}f\big)(a)\cdot
    \|v\|_{L^{p(\cdot)}([0,a])}<\infty.
\end{equation*}
Let us use the following representation for  $M^{(3)}_{\alpha}f(x)$:
\begin{eqnarray*}
% \nonumber to remove numbering (before each equation)
  \big(M^{(3)}_{\alpha}f\big)(x) &=& \chi_{(a,2a]}(x)M^{(\mathbb{R}_+)}_{\alpha}\big(f\cdot\chi_{[0,a]}\big)(x)+
  \chi_{(2a,\infty)}(x)M^{(\mathbb{R}_+)}_{\alpha}\big(f\cdot\chi_{[0,a]}\big)(x).\\
   &=:&\big(\overline{M}^{(3)}_{\alpha}f\big)(x)+\big(\widetilde{M}^{(3)}_{\alpha}f\big)(x).
\end{eqnarray*}
It is easy to check that for $x\in(a,2a]$,
\begin{equation*}
    \big(\overline{M}^{(3)}_{\alpha}f\big)(x)\leq\sup\limits_{h>a-x}\frac{1}{(a-x+h)^{1-\alpha}}
    \int\limits_{x-h}^{a}|f(y)|dy\leq\big(M^{(\mathbb{R}_+)}_{\alpha}f\big)(a).
\end{equation*}
Consequently,
\begin{equation*}
    \|v\overline{M}^{(3)}_{\alpha}f\|_{L^{p(\cdot)}(\mathbb{R}_+)}\leq\|f\|_{L^{p_c}\big((a,2a]\big)}\big
    (M^{(\mathbb{R}_+)}_{\alpha}f\big)(a)<\infty,
\end{equation*}
because $v^{p(\cdot)}(\cdot)$ is locally integrable on $\mathbb{R}_+$. Further we have that for $x>2a$,
\begin{equation*}
    \big(\widetilde{M}^{(3)}_{\alpha}f\big)(x)\leq\frac{1}{(x-a)^{1-\alpha}}\int\limits_{0}^{a}|f(y)| dy.
\end{equation*}
Hence, by using H\"older's inequality  in $L^{p(\cdot)}$ spaces, we find that
\begin{eqnarray*}
% \nonumber to remove numbering (before each equation)
  \bigg\|v\widetilde{M}^{(3)}_{\alpha}f\bigg\|_{L^{p(\cdot)}(\mathbb{R}_+)} &\leq& \bigg\|\frac{v(x)}{(x-a)^{1-\alpha}}\bigg\|_{L^{p_c}\big((2a,\infty)\big)}\bigg( \int\limits_{0}^{a}|f(y)|dy\bigg)   \\
   &\leq&\bigg\|\frac{v(x)}{(x-a)^{1-\alpha}}\bigg\|_{L^{p_c}\big((2a,\infty)\big)} \\
   &&\big\|fw\big\|_{L^{p(\cdot)}\big((0,a]\big)}\big\|w^{-1}\big\|_{L^{p'(\cdot)}\big((0,a]\big)}\\
   &=&I_1\cdot I_2\cdot I_3.
\end{eqnarray*}
Since $I_2<\infty$ and $I_3<\infty$, we need to show that $I_1<\infty$. This follows from the fact that condition $(1.3)$ yields
\begin{equation*}
    \big\|v\overline{M}_{\alpha}\big(w^{-(p_c)'}\chi_{I}\big)\big\|_{L^{p_c}\big((2a,\infty)\big)}\leq
    \big\|w^{1-(p_c)'}(\cdot)\chi_{I}(\cdot)\big\|_{L^{p_c}\big((2a,\infty)\big)},\;\;  I\subset (2a, \infty),    \eqno{(1.4)}
\end{equation*}
where $\overline{M}_{\alpha}$ is the maximal operator defined on $(2a,\infty)$ as follows:
\begin{equation*}
    \big( \overline{M}_{\alpha}f\big) (x)=\sup\limits_{h>0}\frac{1}{h^{1-\alpha}}\int\limits_{(2a,\infty)\cap(x-h,x+h)}|f(y)|dy.
\end{equation*}
Using the result by E. Sawyer  see \cite{Saw} (see also \cite{GaRu}, Ch. 4)  for Lebesgue spaces with constant parameter, we see that
$(1.4)$ implies the inequality
\begin{equation*}
    \big\|v\overline{M}_{\alpha}f\big\|_{L^{p_c}\big((2a,\infty)\big)}\leq c\big\|fw\big\|_{L^{p_c}\big((2a,\infty)\big)}.
\end{equation*}
Since
\begin{equation*}
    \overline{M}_{\alpha}f(x)\geq\frac{1}{(x-a)^{1-\alpha}}\int\limits_{2a}^{x}|f(y)|dy \ \ \ \hbox{for}\ \ \  x>2a,
\end{equation*}
we have that for the Hardy operator
 \begin{equation*}
    \big(H_af\big)(x)=\int\limits_{2a}^{x}f(t)dt, \ \ \ \ \ \ x>2a,
 \end{equation*}
 the two-weight inequality
 \begin{equation*}
    \big\|v(x)(x-a)^{\alpha-1}H_af\big\|_{L^{p_c}\big((2a,\infty)\big)}
    \leq\big\|wf\big\|_{L^{p_c}\big((2a,\infty)\big)} \eqno{(1.5)}
 \end{equation*}
 holds. Let us recall that (see e.g. \cite{Maz}, Section 1.3) necessary condition for $(1.5)$ is that
 \begin{equation*}
    \sup\limits_{t>2a}\bigg(\int\limits_{t}^{\infty}
    \bigg[\frac{v(x)}{(x-a)^{1-\alpha}}\bigg]^{p_c}dx\bigg)^{\frac{1}{p_c}}\bigg(\int\limits_{2a}^{t}
    w^{1-(p_c)'}(x)dx\bigg)^{\frac{1}{(p_c)'}}<\infty.
 \end{equation*}
 Hence,
 \begin{equation*}
  \int\limits_{2a}^{\infty}
    \bigg[\frac{v(x)}{(x-a)^{1-\alpha}}\bigg]^{p_c}dx  =  \int\limits_{2a}^{3a}
  (\cdots) +  \int\limits_{3a}^{\infty}
  (\cdots)  $$
  $$ \leq a^{\alpha-1} \int\limits_{2a}^{3a} \big(v(y)\big)^{p_c}+ \int\limits_{3a}^{\infty}
  \bigg[\frac{v(x)}{(x-a)^{1-\alpha}}\bigg]^{p_c}dx <\infty.
 \end{equation*}
 It remains to estimate $I:=\|vM^{(4)}_{\alpha}f\|_{L^{p(\cdot)}(\mathbb{R}_+)}$. But $I<\infty$ because of the two-weight result by E. Sawyer \cite{Saw} (see also \cite{GaRu}, Ch.4)  for the maximal operator defined on $(a,\infty)$ in Lebesgue spaces with constant
 exponent. Sufficiency is proved.

 {\em Necessity} follows easily by taking the test functions $f(\cdot)= \chi_{I}(\cdot) w^{-p'(\cdot)}(\cdot)$ in the two--weight inequality.
 \end{proof}

\vskip+0.2cm
The next statement follows in the same way as the previous one; therefore we omit the proof.
\vskip+0.2cm
{\bf Theorem 1.3.} {\em Let $0\leq \alpha<1$, $1<p_-\leq p\leq p_+<\infty$,  and let $p\in LH(\mathbb{R})$. Suppose that there is a positive number $a$ such that
$w^{-p'(\cdot)}(\cdot)\in DC([-a,a])$ and $p\equiv p_c\equiv $const outside $[-a,a]$. Then the inequity
\begin{equation*}
    \|v M^{(\mathbb{R})}_{\alpha}f\|_{L^{p(\cdot)}(\mathbb{R})}\leq\|wf\|_{L^{p(\cdot)}(\mathbb{R})},
\end{equation*}
holds  if and only if there is a positive constant $b$ such that for all bounded intervals $I\subset\mathbb{R} $,
\begin{equation*}
    \|v M^{(\mathbb{R})}_{\alpha}(w^{-p'(\cdot)}\chi_{I})\|_{L^{p(\cdot)}(\mathbb{R})}\leq c\|w^{1-p'(\cdot)}\|_{L^{p(\cdot)}(I)}<\infty.
\end{equation*}}

\section{Integral operators on ${\Bbb{R}}_+$}

 In this section we derive  two--weight criteria of  other type for the operators

$$ ({\cal{H}}f)(x)= \text{(p.v.)} \; \int\limits_0^{\infty} \frac{f(t)}{x-t}dt, \;\; x\in {\Bbb{R}}_+, $$

$$ ({\mathcal{M}}f)(x) = \sup_{I\ni x} \frac{1}{|I|} \int\limits_I |f(t)| dt, \;\; x\in {\Bbb{R}}_+, $$
provided that weights are monotonic,  where the supremum is taken over all finite intervals $I\subset
{\Bbb{R}}_+$ containing $x$.

In this section we shall use the notation
$$g_-:= g_-({\Bbb{R}}_+); \;\;\; g_+:= g_+({\Bbb{R}}_+), $$
for a measurable function $g: {\Bbb{R}}_+\to {\Bbb{R}}_+$.

\vskip 0.3cm

First we present the following statement regarding the weighted
Hardy transform
$$ (H_{v,w}f)(x)= v(x) \int\limits_0^x f(t)w(t) dt $$
and its dual
$$ (H'_{v,w}f)(x)= v(x) \int\limits_x^{\infty} f(t)w(t) dt$$
defined on ${\Bbb{R}}_+$.

\vskip+0.2cm

{\bf Theorem A.} Let $1<p_-\leq p(x) \leq q(x) \leq q_- <\infty$
and let $p,q\in LH({\Bbb{R}}_+)$. Suppose that $p=p_c\equiv
const$, $q=q_c\equiv const$ outside some interval $(0,a)$. Then

\rm{(i)} the operator $H_{v,w}$ is bounded from
$L^{p(\cdot)}({\Bbb{R}}_+)$ to $L^{q(\cdot)}({\Bbb{R}}_+)$ if and
only if
$$ D:= \sup_{t>0}D(t):= \sup_{t>0}\| v \|_{L^{q(\cdot)}\big((t,\infty)\big)} \| w \|_{L^{p'(\cdot)}\big( (0,t) \big)}< \infty;$$

\rm{(ii)} the operator $H'_{v,w}$ is bounded from
$L^{p(\cdot)}({\Bbb{R}}_+)$ $L^{q(\cdot)}({\Bbb{R}}_+)$ if and
only if
$$ D':= \sup_{t>0}D'(t):= \sup_{t>0}\| v \|_{L^{q(\cdot)}\big( (0,t) \big)} \| w \|_{L^{p'(\cdot)}\big(  (t, \infty) \big)}< \infty.$$
\vskip+0.1cm

{\em Proof.} We prove part (i). Part (ii) follows from the duality arguments.
Let $\|f\|_{L^{q(\cdot)}({\Bbb{R}}_+)}\leq 1$. We represent $H_{v, w}f$ as follows:

$$H_{v,w}f(x)= \chi_{[0,a]}v(x) \int\limits_0^x f(t) w(t) dt +
\chi_{(a,\infty)}v(x) \int\limits_0^x f(t) w(t) dt :=
H^{(1)}_{v,w}f(x)+H^{(2)}_{v,w}f(x).$$

Observe that the condition $D<\infty$ implies that

$$ D^{(a)}:= \sup_{0<t<a} \| v \|_{L^{q(\cdot)}\big((t,a)\big)} \| w \|_{L^{p'(\cdot)}\big( (0,t) \big)}< \infty.$$
Consequently (see \cite{Kop}),

$$\|H^{(1)}_{v,w}f  \|_{L^{q(\cdot)}({\mathbb{R}})} \leq c \| f  \|_{L^{p(\cdot)}([0,a])} \leq c. $$

It remains to estimate $\|H^{(2)}_{v,w}f
\|_{L^{q(\cdot)}({\Bbb{R}}_+)}$. Let $\|g
\|_{L^{q'(\cdot)}({\Bbb{R}}_+)}\leq 1$. We have that
$$\int\limits_0^{\infty} (H^{(2)}_{v,w}f)(x) g(x) dx = \int\limits_a^{\infty}
(H^{(2)}_{v,w}f)(x) g(x)   dx $$
$$ \leq \int\limits_a^{\infty}  v(x) \bigg( \int\limits_a^x f(t) w(t) dt \bigg)  g(x)   dx + \bigg(\int\limits_a^{\infty}  v(x) g(x) dx\bigg)  \bigg( \int\limits_0^a f(t) w(t) dt \bigg):= S_1 +S_2. $$

We can now apply the boundedness of the Hardy transform
$T^{(a)}_{v,w}f(x)= v(x) \int\limits_a^x f(t) w(t) dt$ from
$L^{p_c}([a,\infty))$ to $L^{q_c}([a,\infty))$  (see e.g. \cite{Maz}, Section 1.3) because
$$ \sup_{t>a} \| v \|_{L^{q_c}\big((t,\infty)\big)} \| w \|_{L^{(p_c)'}\big( (a,t) \big)}\leq D<\infty.$$
Consequently, by this fact and H\"older's inequality we derive that
$$ S_1 \leq \| T^{(a)}_{v,w}f \|_{L^{q_c}([a,\infty))} \|g\|_{L^{q_c}([a,\infty))} \leq c \| f \|_{L^{p(\cdot)}({\Bbb{R}}_+)} \leq C. $$

Applying H\"older's inequality for $L^{p(\cdot)}$ spaces we find that
$$  S_2 \leq \bigg(\int\limits_a^{\infty}  v(x) g(x) dx\bigg)  \| f\|_{L^{p(\cdot)}([0,a])} \|w\|_{L^{p'(\cdot)}([0,a])} \leq C. $$

%We construct the sequence $\{ t_k \}$ so that $2^k = \int\limits_{0}^{t_k} f(x) dx $. Then it is easy to see that $2^k = \int\limits_{t_k}^{t_{k+1}} %f(x) dx $. Let $g\in L^{q'(\cdot)}({\Bbb{R}}_+)$ with $\|g\|_{L^{q'(\cdot)}({\Bbb{R}}_+)}\leq 1$. Using these argument, the notation $E_k:= [t_k, %t_{k+1})$ and Lemma 3.1 we have that
%$$ \int\limits_0^{\infty} (H_{v,w}f)(x) g(x) dx \leq \sum_{k}\Big(\int_{t_k}^{t_{k+1}} g(x) v(x)  dx \Big) \Big(\int\limits_0^{t_{k+1}} f(x) w(x) %dx\Big) $$
%$$ = c \sum_{k}\Big(\int_{E_k} g(x) v(x)  dx \Big) \Big(\int\limits_{E_{k-1}} f(x) w(x) dx\Big)$$
%$$ \leq  \sum_{k} \| f \chi_{E_{k-1}} \|_{L^{p(\cdot)}( {\Bbb{R}}_+)} \| g \chi_{E_k}
%\|_{L^{q'(\cdot)}( {\Bbb{R}}_+)} \leq  c \| f
%\|_{L^{p(\cdot)}( {\Bbb{R}}_+)} \| g \|_{L^{q'(\cdot)}( {\Bbb{R}}_+)}\leq 1.$$

{\em Necessity} follows by the standard way choosing the appropriate
test functions. $\Box$

\vskip+0.1cm

{\bf Theorem B (\cite{EdKoMe2}).} {\em   $1<p_-\leq p_+<\infty$.
Suppose that $p\in LH({\Bbb{R}}_+)$ and that $p=p_c=const$ outside
some interval. Then the inequality

$$ \| v Tf
\|_{L^{p(\cdot)}({\Bbb{R}}_+)}  \leq c \| w
f\|_{L^{p(\cdot)}({\Bbb{R}}_+)}, \eqno{(2.1)}
$$
where  $T$ is ${\cal{M}}$ or ${\cal{H}}$, holds if

\rm{(i)} $\;\;\;\; H_{\overline{v}, \widetilde{w}}$ is bounded in
$L^{p(\cdot)}({\Bbb{R}})$, where $\overline{v}(x):=
\frac{v(x)}{x}$, $\widetilde{w}(x) := \frac{1}{w(x)}$;

\rm{(ii)} $\;\;\;\; H'_{v, \widetilde{w}_1}$ is bounded in
$L^{p(\cdot)}({\Bbb{R}})$, where $\widetilde{w}_1(x):=
\frac{1}{w(x) x}$;

\rm{(iii)}  $$ v_+([x/4, 4 x]) \leq c w(x) \;\; \text{a.e. or}\;\;
v(x)\leq c w_-([x/4, 4x])\; \text{a.e.} \eqno{(2.2)}
$$}

\vskip+0.1cm

Theorems A and B imply the following statement:

\vskip+0.1cm

{\bf Theorem 2.1.} {\em  Let $1<p_-\leq p_+ <\infty$ and let $p\in
LH({\Bbb{R}}_+)$. Suppose that $p=p_c\equiv \;const$ outside some
interval $[0,a]$. Suppose also that $v$ and $w$ are weights on
${\Bbb{R}}_+$. Then the inequality $(2.1)$, where  $T$ is
${\cal{M}}$ or ${\cal{H}}$, holds if

\rm{(i)} $$ E_1:= \sup_{t>0}E_1(t):= \sup_{t>0}\| v(x)x^{-1}
 \|_{L^{p(x)}\big((t,\infty)\big)} \| w^{-1} \|_{L^{p'(\cdot)}\big( (0,t) \big)}<
 \infty;\eqno{(2.3)}$$

\rm{(ii)} $$ E_2:= \sup_{t>0}E_2(t):= \sup_{t>0}\| v
\|_{L^{p(\cdot)}\big((0,t)\big)} \| w^{-1}(x)x^{-1}
\|_{L^{p'(x)}\big( (0,t) \big)}< \infty; \eqno{(2.4)}$$

\rm{(iii)} condition $(2.2)$ is satisfied.}

\vskip+0.1cm

Now we prove the next statement.

\vskip+0.1cm

{\bf Theorem 2.2.} {\em  Let $1<p_-\leq p_+ <\infty$ and let $p\in
LH({\Bbb{R}}_+)$. Suppose that $p=p_c\equiv \;const$ outside some
interval $[0,a]$. Suppose also that $v$ and $w$ are positive
increasing functions on ${\Bbb{R}}_+$. Then inequality $(2.1)$,
where  $T$ is ${\cal{M}}$ or ${\cal{H}}$, holds if and only if
$(2.3)$ is satisfied. }

\vskip+0.1cm

{\em Proof.}  {\em Sufficiency.} Taking Theorem 2.1 into account it
is enough to see that condition (2.3) implies conditions (2.4) and
(2.2). For (2.2) we will show that there is a positive constant
$c$ such that for all $t>0$  inequality $$ v(4t) \leq c w(t),
\;\;\; t>0. \eqno{(2.5)}$$
holds. Indeed,  inequality (1.1) with respect to the Lebesgue measure $d\mu(x)= dx$ and the exponent $r=p'$ which belongs to $LH([0,a])$, for small $t$, yields that

 $$  E_{1}(t) \geq \|\chi_{[t, 4t]}(\cdot)|\cdot|^{-1}\|_{L^{p(\cdot)}_{v(\cdot)}(\mathbb{R}_+)}
\parallel \chi_{[0, t/4]}(\cdot)w^{-1}(\cdot)\parallel _{L^{p'(\cdot)}(\mathbb{R}_+)}$$

$$ \geq c\frac{v(t)}{t} t^{\frac{1}{p_-([t,4t])}}
w^{-1}(t/4)t^{\frac{1}{(p')_-([0,t/4])}} \geq
  c\frac{v(t)}{w(t/4)}t^{-1} t^{\frac{1}{p_-([0,4t])}} t^{\frac{1}{(p')_-([0,t/4])}}= c\frac{v(t)}{w(t/4)}. $$
Further, for large $t$, we have that

$$ E_{1}(t) \geq \| v(x) x^{-1} \chi_{(t, 2t)}(x)\|_{L^{p_c}(\mathbb{R}_+)} \|
 \chi_{[t/8 , t/4]}(\cdot)w^{-1}(\cdot)\|_{L^{p'_c}(\mathbb{R}_+)} \geq c\frac{v(t)}{w(t/4)}t^{-1} t^{\frac{1}{p_c}}t^{\frac{1}{(p_c)'}}=c\frac{v(t)}{w(t/4)} $$
Thus, condition  $(2.2)$ is satisfied.

Taking into account the fact that $v$ and $w$ are increasing and
inequality (2.5) we can easily conclude  that  condition (2.4) is
satisfied.

{\em Necessity}. First observe that inequality (2.1) implies that
$\|w^{-1}\|_{L^{p'(\cdot)}(0,t)} <\infty $ for all $t>0$.

Let $T={\cal{M}}$. Then using the obvious inequality
$$ {\cal{M}}f(x) \geq \frac{c}{x} \int_0^x f(t) dt, \;\;\; x>0, $$
and taking into account Theorem A we have necessity for
${\cal{M}}$. Let now $T= {\cal{H}}$. We take $f \geq 0$ so that
$\| f\|_{L^{p(\cdot)}_w({\Bbb{R}}_+)}\leq 1$.  Then we have that
$$
\|v {\cal{H}}f\|_{L^{q(\cdot)}({\Bbb{R}}_+)} \leq C. \eqno{(2.6)}
$$
Obviously, (2.6) yields that
$$ C \geq \|v {\cal{H}}f\|_{L^{q(\cdot)}({\Bbb{R}}_+)} \geq \| \chi_{(t, \infty)}(\cdot) v {\cal{H}}f \|_{L^{p(\cdot)}({\Bbb{R}}_+)}.$$
If $f$ has support on $(0,t)$, $t>0$,  then this inequality
implies  that

$$  C \geq  \bigg\| \chi_{(t, \infty)}(\cdot) v(\cdot) \bigg( \int_{0}^t  \frac{f(y)}{\cdot-y}dy\bigg) \bigg\|_{L^{p(\cdot)}({\Bbb{R}}_+)} \geq c \bigg\| \chi_{(t, \infty)}(x) v(x)x^{-1}  \bigg\|_{L^{p(\cdot)}({\Bbb{R}}_+)}\bigg( \int_{0}^t  f(y)dy\bigg). $$

%Analogously it follows that
%$$ \|v {\cal{H}}f\|_{L^{q(\cdot)}({\Bbb{R}}_+)} \geq c \| \chi_{(-\infty, -t)}(x) v(x)(-x)^{-1}  \bigg\|_{L^{p(\cdot)}({\Bbb{R}}_+)}\bigg( \int_{-t}^t %dy f(y)dy\bigg).$$
%Summarazing these inequalities we conclude that
%$$   C \geq \|v {\cal{H}}f\|_{L^{q(\cdot)}({\Bbb{R}}_+)}\geq c \| \chi_{\{ |\cdot|>t\}} (x) v(x)x^{-1}  \bigg\|_{L^{p(\cdot)}({\Bbb{R}}_+)}\bigg( %\int_{-t}^t dy f(y)dy\bigg)$$.
By taking now supremum with respect to $f$ and using the
inequality
$$\|g\|_{L^{p(\cdot)}}\leq \sup_{\|h\|_{L^{p'(\cdot)}}\leq 1}\bigg| \int g h\bigg|,  $$
(see e.g. \cite{Sa1}) we have necessity. $\Box$.

\vskip+0.5cm

\section*{Acknowledgements}
The authors were partially supported by the
Georgian National Science Foundation Grant (project numbers: No.
GNSF/ST09/23/3-100 and No. GNSF/ST07/3-169).    \vskip+1cm

\vskip+0.5cm

Authors' Addresses:

\

V. Kokilashvili: A. Razmadze Mathematical Institute, 1. M. Aleksidze Str., 0193 Tbilisi, Georgia and Faculty of Exact and Natural Sciences, I. Javakhishvili Tbilisi State University 2,
University St., Tbilisi 0143 Georgia \\ e-mail: kokil@rmi.acnet.ge.\\

\

A. Meskhi: A. Razmadze Mathematical Institute, 1. M. Aleksidze Str., 0193 Tbilisi, Georgia and Department of Mathematics,  Faculty of Informatics and Control Systems, Georgian Technical University, 77, Kostava St., Tbilisi, Georgia.\\e-mail: meskhi@rmi.acnet.ge\\

\end{document}